\documentclass[12pt]{article}
\begin{document}

\title{On the Asymptotic Formula for the Number of Plane
Partitions of Positive Integers}

\author{\bf  Ljuben Mutafchiev
\\ American University in Bulgaria
and
\\ Institute of Mathematics and Informatics of the Bulgarian Academy
of Sciences
\\ \tt {ljuben@aubg.bg}
\\ and
\\ \bf Emil Kamenov
\\ Sofia University
\\ \tt {kamenov@fmi.uni-sofia.bg} }

\date{ }

\renewcommand{\baselinestretch}{1.3}

\maketitle

A plane partition $\omega$ of the positive integer $n$ is an array
of non-negative integers
\begin{equation}
      \begin{array}{lllc}
      \omega_{1,1} & \omega_{1,2} & \omega_{1,3} & ... \\
      \nonumber
      \omega_{2,1} & \omega_{2,2} & \omega_{2,3} & ... \\
      ... & ... & ... & ...\\
      \nonumber
      \end{array}
\end{equation}
 for which
 $$
 \sum_{i,j}\omega_{i,j}=n
 $$
 and the rows and columns are arranged in decreasing order:
 $$
 \omega_{i,j}\ge\omega_{i+1,j}, \omega_{i,j}\ge\omega_{i,j+1}
 $$
 for all $i,j\ge 1$. The non-zero entries $\omega_{i,j}>0$ are
 called parts of $\omega$. Sometimes, for the sake of
 brevity, the zeroes in the array (1) are deleted. For instance,
 the abbreviation
 $$
 \begin{array}{lll}
 3 & 2 & 1 \\
 1 & 1 & \\
 \end{array}
 $$
 is assumed to present a plane partition of $n=8$ having $2$
 rows and $5$ parts.

 It seems that MacMahon was the first who introduced the idea of a
 plane partition; see [5]. He deals with the general problem of
 such partitions enumerating them by the size of each part, number
 of rows and number of columns. These problems have been
 subsequently reconsidered by other authors who have developed
 methods, entirely different from those of MacMahon. For important
 references and more details in this direction we refer the reader
 to the monograph of Andrews [1; Chap. 11] and to the survey paper
 of Stanley [8; Chap. V].

 Let $q(n)$ denote the total number of plane partitions of the
 integer $n\ge 1$. It turns out that
 \begin {equation}\label{genf}
 Q(x)=1+\sum_{n=1}^\infty q(n)x^n=\prod_{j=1}^\infty(1-x^j)^{-j}
 \end{equation}
 (see [1; Corollary 11.3] or [8; Corollary 18.2]). The asymptotic
 of $q(n)$ has been obtained by Wright [10]. It is given by the
 following formula:
\begin{eqnarray}\label{qofn}
& & q(n)=\frac{[\zeta(3)]^{7/36}}{2^{11/36}\pi^{1/2}} n^{-25/36}
\exp{\{3[\zeta(3)]^{1/3}(n/2)^{2/3}+2c\}} \nonumber \\ & & \times
[\gamma_0+\sum_{h=1}^r\gamma_h n^{-2h/3} +O(n^{-2(r+1)/3})],
\end{eqnarray}
where
 $$
 \zeta(s)=\sum_{j=1}^\infty j^{-s}
 $$
is the Riemann's zeta function,
\begin{equation}
c=\int_0^\infty\frac{y\log{y}}{e^{2\pi y}-1} dy
 \end{equation}
and $\gamma_h$ are certain constants with explicitly given values.
(\ref{qofn}) implies the asymptotic equivalence
\begin{equation}\label{qofn2}
q(n)\sim\frac{\gamma_0[\zeta(3)]^{7/36}}{2^{11/36}\pi^{1/2}}
n^{-25/36}\exp{\{3[\zeta(3)]^{1/3}(n/2)^{2/3}+2c\}}, n\to\infty.
\end{equation}
Many authors refer to (\ref{qofn}) and (\ref{qofn2}) in their work
on this subject. Wright [10; p. 179] states in his theorem that
$\gamma_0=1$, however, at the the end of its proof his argument
shows that the true value of $\gamma_0$ is $\gamma_0=3^{-1/2}$. In
fact, the substitution $v=1+y$ [10; p. 188] that he makes in the
integrals denoted by $a_{2m}$ generates the factor
\begin{equation}\label{eq6}
(3+2y)^{-(m+1/2)}, m=0,1,...
\end{equation}
in their integrands. Then, Wright shows how $a_{2m}$ determine the
coefficients $\gamma_h$ in (\ref{qofn}). He proves that
$a_{2m}=(-1)^m b_{sm}/2\pi$, where $b_{sm}$ is the coefficient of
$y^{2m}$ in the expansion of
 $$
 (1+y)^{2s+2m+13/12}(3+2y)^{-m-1/2}
 $$
(here $s$ is a non-negative integer; see also the definition of
the constants $b_{sm}$ on p. 179 [10; formula (2.15)]). It turns
out that if $h=0$, one has to set in (\ref{eq6}) $m=0$, and thus
by [10; formula (2.22)], $\gamma_0$ equals the constant term in
the power series expansion of
 $$
 (1+y)^{2s+13/12}(3+2y)^{-1/2}.
 $$
Hence $\gamma_0=3^{-1/2}$.

The asymptotic equivalence (\ref{qofn2}) can be also obtained
after a direct application of a theorem due to Meinardus [6] (see
also [1; Chap. 6]) who has obtained the asymptotic behavior of the
coefficients in the power series expansion of infinite products of
the following form:
 $$
 \prod_{j=1}^\infty(1-x^j)^{-a_j},
 $$
where $\{a_j\}_{j\ge 1}$ is given sequence of non-negative
numbers. He introduced a scheme of assumptions on $\{a_j\}_{j\ge
1}$, which are satisfied by the generating function (\ref{genf})
of the numbers $q(n)$. Further generalizations of Meinardus'
result are given in [3].

The object of this work is to show that Meinardus' theorem implies
asymptotic formula (\ref{qofn2}) with $\gamma_0=3^{-1/2}$. This
fact was also briefly mentioned in [4], where the limiting
distribution of the trace of a random plane partition was studied.
A sketch of the proof of (\ref{qofn2}) based on the fact that
$Q(x)$ is an admissible function in the sense of Hayman [2] was
presented there.

We start with the statement of Meinardus theorem [6] (see also [1;
Chap. 6]) . Let

$$
 \prod_{j=1}^\infty \left( 1 - x^j \right)^{-a_j} = 1
+ \sum_{n=1}^\infty r(n) x^n,
$$
 where $x = e^{-\tau}$, $\Re e\
\tau
>0$. Consider also the
associated Dirichlet series

\begin{equation}\label{dofs}
D(s) = \sum_{j=1}^\infty \frac{a_j}{j^s}, \hspace{.6cm} s=\sigma +
it.
\end{equation}
Meinardus assumes that the following four conditions hold.

{\bf (i)} $D(s)$ converges in the half-plane  $\sigma > \alpha >
0$.

{\bf (ii)} $D(s)$ can be analytically continued in the region $\Re
e\ s\geq -C_0$ $(0 < C_0 < 1)$, where $D(s)$ is analytic except
for a simple pole at $s=\alpha$ with residue $A$.

{\bf (iii)} There exists a constant $C_1>0$ such that $D(s)=
O(|t|^{C_1})$ uniformly in $\sigma \geq - C_0$ as $|t| \rightarrow
\infty$.

{\bf (iv)} Let
\begin{equation}\label{gofv}
g(v)=\sum\limits_{j=1}^\infty a_j x^j, \hspace{.6cm} x= e^{-v},
\end{equation}
where $v=y + 2\pi i w$ and $y$ and $w$ are real numbers. For
$|\arg v|
> \pi/4$ and $|w|\leq 1/2$ and for sufficiently small $y$, $g(v)$
satisfies
\begin{equation}\label{eq9}
\Re e\ g(v)-g(y) \leq -C_2 y^{-\varepsilon},
\end{equation}
where $\varepsilon>0$ is an arbitrary number, and $C_2
>0$ is suitably chosen and may depend on $\varepsilon$.

{ \bf Theorem} [6] {\sl If conditions (i) - (iv) hold, then $$
r(n) = C n^K \exp\left\{n^{\frac{\alpha}{\alpha+1}} \left(1 +
\frac{1}{\alpha} \right)  [A \Gamma(\alpha+1) \zeta(\alpha+1)
]^{\frac{1}{\alpha+1}} \right\} (1 + O(n^{-K_1})), $$ where
$\Gamma(\alpha)$ is Euler's gamma function,
\begin{eqnarray}\label{C}
&&C = e^{D'(0)} [2 \pi (1 + \alpha)]^{-1/2} [A \Gamma(\alpha+1)
\zeta(\alpha+1)]^{\frac{1 -2 D(0)} {2 + 2 \alpha}},\hspace{3.cm}\\
\label{k}
 &&K =\frac{D(0)  -1 - \alpha/2}{ \alpha + 1},\\ \label{k1}
 &&K_1 = \frac{\alpha}{\alpha +1}\min \left( \frac{C_0}{\alpha} - \frac{\delta}{4}, \frac{1}{2}
 - \delta \right)
 \end{eqnarray}
and $\delta$ is an arbitrary positive number.

}

In the case of plane partitions we have $\{a_j\}_{j\ge 1} =
\{j\}_{j\ge 1}$ and thus by (\ref{dofs}), $D(s) = \zeta(s-1)$.
Thus, (i) - (iii) follow from well known properties of the zeta
function (see e.g. [9; Chap. 13]). Moreover, in (ii) $C_0$ is an
arbitrary constant within range $(0,1)$, since the zeta function
has analytical continuation in whole complex plane. It is easily
seen that $\alpha = 2$, $A = 1$. To verify condition (iv) notice
first that for $\{a_j\}_{j\ge 1}=\{j\}_{j\ge 1}$ the function
$g(v)$ defined by (\ref{gofv}) becomes
 $$
 g(v)=\sum_{j=1}^\infty j\exp{(-jy+2\pi ijw)}.
 $$
Hence
 $$
 \Re e g(v) - g(y)=\sum_{j=1}^\infty je^{-jy}[\cos{(2\pi
jw)} -1]
 =\sum_{k=1}^\infty(-1)^k\frac{(2\pi w)^{2k}}{(2k)!}
 \sum_{j=1}^\infty j^{2k+1}e^{-jy}.
 $$
We simplify the inner sum here replacing its terms by the
following derivatives:
 $$
 j^{2k+1}e^{-jy}=(-1)^{2k+1}\frac{d^{2k+1}}{dy^{2k+1}}(e^{-jy}).
 $$
In this way we get
 \begin{eqnarray}\label{eq13}
 & & \Re e g(v) - g(y)
 =\sum_{k=1}^\infty (-1)^{-(k+1)}\frac{(2\pi w)^{2k}}{(2k)!}
 \frac{d^{2k+1}}{dy^{2k+1}}(\sum_{j=1}^\infty e^{-jy}) \nonumber
 \\
 & & =\sum_{k=1}^\infty (-1)^{-(k+1)}\frac{(2\pi w)^{2k}}{(2k)!}
 \frac{d^{2k+1}}{dy^{2k+1}}\left(\frac{1}{e^y-1}\right).
 \end{eqnarray}
 Further, to calculate the $(2k+1)$th derivative of $1/(e^y-1)$
 we apply Faa di Bruno's formula [7; Section 2.8]. Using also an
 identity for the Stirling numbers of the second kind $\sigma_n^{(m)}$
 [7; Section 4.5],
 we find that:
 \begin{eqnarray}
 & & \frac{d^{2k+1}}{dy^{2k+1}}\left(\frac{1}{e^y-1}\right)
 =\sum_{m=0}^{2k+1}\frac{(-1)^m m!}{(e^y-1)^{m+1}} \sum_b
 \frac{(2k+1)!e^{my}}{(1!)^{b_1}b_1!(2!)^{b_2}...[(2k+1)!]^{b_{2k+1}}}
 \nonumber \\
 & & =\sum_{m=0}^{2k+1}\frac{(-1)^m m!e^{my}}{(e^y-1)^{m+1}}
 \sigma_{2k+1}^{(m)} =(e^y-1)^{-(2k+2)} \sum_{m=0}^{2k+1}(-1)^m
 m!e^{my}(e^y-1)^{2k+1-m}\sigma_{2k+1}^{(m)},\nonumber
 \end{eqnarray}
where $\sum_b$ means that the summation is over all non-negative
integers $b_j$ that satisfy $n=\sum_j jb_j, m=\sum_j b_j$. It is
clear that all terms in the last sum are close to zero if $y$ is
sufficiently small except the last one. It equals
$(-1)^{2k+1}(2k+1)!\sigma_{2k+1}^{(2k+1)}=(-1)^{2k+1}(2k+1)!$.
Therefore, after simple algebraic manipulations, we can rewrite
(\ref{eq13}) in the following form:
 \begin{eqnarray}\label{eq14}
 & & \Re e g(v)-g(y)=\sum_{k=1}^\infty(-1)^k(2\pi w)^{2k}
 (e^y-1)^{-(2k+2)}[2k+1+\psi_k(y)] \nonumber \\
 & & =-\frac{(2\pi w)^2}{(e^y-1)^4}[3+\psi_1(y)]
 +\left\{\frac{(2\pi w)^4}{(e^y-1)^6}[5+\psi_2(y)] -\frac{(2\pi
 w)^6}{(e^y-1)^8}[7+\psi_3(y)]\right\}+...,
\end{eqnarray}
where
 $$
 \psi_k(y)=\frac{(-1)^{2k+1}}{(2k)!}\sum_{m=0}^{2k}(-1)^m
 m!e^{my}(e^y-1)^{2k+1-m}\sigma_{2k+1}^{(m)}=O(y)
 $$
 as $y\to0^{+}$. For sufficiently small $y$ we also have
 $e^y-1=y+O(y^2)$.
 Moreover, the requirement $\mid\arg\tau\mid\ge\pi/4$ implies
 that $(2\pi w)^{2k}\ge y^{2k}$. Using these arguments,
 we conclude
 that the first term in (14) is $-(3/y^2)[1+O(y)]$ as $y\to
 0^{+}$. For the other terms in the curly brackets of (\ref{eq14}) we get
 the estimate
 \begin{equation}\label{eq15}
 \frac{(2\pi w)^{4l}}{y^{4l+2}}\left\{4l+1-\left(\frac{2\pi
 w}{y}\right)^2[4l+3+O(y)]\right\}\le\frac{(2\pi
 w)^{4l}}{y^{4l+2}}[-2+O(y)]\le-2y^{-2}+O(y^{-1}),
 \end{equation}
 where $l=1,2,...$. Suppose now that $\epsilon\in(0,2]$.
 Estimate (15) implies that one can find certain constant $C_2$
 such that
 $$
 y^{\epsilon}[\Re e g(v)-g(y)]\le-C_2<0.
 $$
 If $\epsilon>2$, then for any $\eta>0$ and sufficiently small $y$
 we have $0<y^{\epsilon-2}<\eta$. Therefore for such $y$'s one
 can define a sequence of enough large positive integers $\{M_y\}$
 satisfying
 the inequality $\eta/2M_y<y^{\epsilon-2}$, or equivalently,
 $-\eta>-2M_y y^{\epsilon-2}$. Now, estimate (\ref{eq15}) shows
 that $-\eta$ exceeds
 the sum of $M_y$ terms of the expansion of $y^{\epsilon}[\Re e
 g(v)-g(y)]$ and this sum
 in turn exceeds the sum of all terms of the expansion.
 Hence, for $\epsilon>2$, we also obtain
 inequality (\ref{eq9}) of condition (iv) with $C_2=\eta$.

 To calculate the constants $D(0)$ and $D^\prime(0)$ we use the
following functional equations (see [9; Chap. 13]):

\begin{equation}\label{zeta}
\zeta(1-z) = 2 \cos\frac{1}{2} \pi z \: (2\pi)^{-z} \Gamma(z)
\:\zeta(z),
\end{equation}

\begin{equation}\label{gama}
\Gamma(z) \:\zeta(z) = \int_0^\infty \frac{w^{z-1}}{e^w-1} dw = (2
\pi)^{z} \int_0^\infty \frac{w^{z-1}}{e^{2 \pi w}-1}  dw.
\end{equation}
Substituting $z=2$ in (\ref{zeta}) we get:

\begin{equation}\label{dof0}
D(0) = \zeta(-1) = - \frac{1}{12}.
\end{equation}
Next, differentiating (\ref{zeta}) and (\ref{gama}) and combining
the expressions of the corresponding derivatives, we obtain
$$
 \zeta'(1-z) = \sin\frac{1}{2} \pi z \: 2^{-z} \pi^{1-z}
\:\Gamma(z)\:\zeta(z) - 2\cos\frac{1}{2} \pi z \: \int_0^\infty
\frac{w^{z-1} \ln w}{e^{2 \pi w}-1}\:  dw,
$$
Substituting again
$z=2$, we find that
\begin{equation}\label{dprime}
D^\prime(0) = \zeta^\prime(-1) = 2 \int_0^\infty \frac{w \ln
w}{e^{2 \pi w}-1}\: dw = 2 c,
\end{equation}
where $c$  is the constant defined by (4).

From (\ref{C}) - (\ref{k1}), (\ref{dof0}) and (\ref{dprime}) we
obtain

$$
 C = 2^{-11/36} (3\pi)^{-\frac{1}{2}}
\left[\zeta(3)\right]^{7/36} e^{2 c},\:\:\:
 K = - \frac{25}{36},\:\:\:
 K_1 = \frac{1}{3} - \varepsilon_1,
$$ where $\varepsilon_1 > 0$. Notice that in (\ref{k1}) we assume
that $\delta$ is sufficiently small and $C_0$ is close to 1. Then,
applying Meinardus theorem, we obtain

$$
 q(n)\sim \frac{[\zeta(3)]^{7/36}}{2^{11/36}3^{1/2}\pi^{1/2}}\
n^{-25/36} \exp \left\{3[\zeta(3)]^{1/3} (n/2)^{2/3} + 2c
\right\},
 $$
 which shows in turn that Wright's formula [10] with
$\gamma_0=3^{-1/2}$ is valid.

We also present some numerical computations of $q(n)$ based on
formula (\ref{qofn2}) and the following recurrence:
 $$
 n\: q(n) = \sum_{k=1}^n q(n-k) \:\beta_2(k).
 $$
 Here $\beta_2(k)$ is the sum of the squares of the positive
divisors of $k$ (see [1; Section 14.6]). We used Maple to get the
results summarized in the table below. The exact values of $q(n)$
are listed in its second column. The third and fourth columns
contain the corresponding values of $q(n)$ computed when
$\gamma_0=1$ and $\gamma_0=3^{-1/2}$, respectively.

\begin{eqnarray}\nonumber
\begin{array}{|r|l|l|l|}
  \hline \phantom{\textrm{\Large A}}n & \phantom{0000}q(n) & \gamma_0=1 & \gamma_0= 3^{-1/2}\\ \hline
  \phantom{\textrm{\Large A}}10      & \phantom{00}500 & \phantom{000}910.69& \phantom{00}525.79 \\ \hline
  \phantom{\textrm{\Large A}}100     & 59\:206 \times 10^{12}  & 103\:709\times 10^{12}   & 59\:876 \times 10^{12} \\ \hline
  \phantom{\textrm{\Large A}}1\:000  & 35\:426 \times 10^{80}  & \phantom{0}61\:507\times 10^{80}    & 35\:511 \times 10^{80} \\ \hline
  \phantom{\textrm{\Large A}}10\:000 & 45\:075 \times 10^{397} & \phantom{0}78\:113\times 10^{397}   & 45\:098 \times 10^{397} \\ \hline
\end{array}
\end{eqnarray}


\begin{thebibliography}{99}


\bibitem{} G. E. Andrews, {\sl The Theory of Partitions, Encyclopedia Math.
Appl. 2}, Addison-Wesley, 1976.
\bibitem{ }W. K. Hayman, {\sl A generalization of Stirling's
formula}, J. Reine Angew. Math. 196(1956), 67-95.
\bibitem{} H.-K Hwang, {\sl Limit theorems for the number of summands
in integer partitions}, J. Comb. Theory Ser. A 96(2001), 89-126.
\bibitem{ } E. Kamenov and L. Mutafchiev, {\sl The limiting
distribution of the trace of a random plane partition}, preprint
available from http://arxiv.org/abs/math.CO/0411377.
\bibitem{ } P. A. MacMahon, {\sl Combinatory Analysis, Vol. 2},
Cambridge Univ. Press, 1916; reprinted by Chelsea, New York, 1960.
\bibitem{ } G. Meinardus, {\sl Assymptotische Aussagen \"{u}ber
Partitionen}, Math. Z. 59(1954), 388-398.
\bibitem{ } J. Riordan, {\sl An Introduction to Combinatorial
Analysis}, J. Wiley, New York, 1958.
\bibitem{ } R. P. Stanley, {\sl Theory and applications of plane
partitions I, II}, Studies Appl. Math. 50(1971), 167-188, 259-279.
\bibitem{ } E. T. Whittaker and G. N. Watson, {\ A Course of
Modern Analysis}, 4th ed., Cambridge Univ. Press, Cambridge, 1927.
\bibitem{ } E. M. Wright, {\sl Asymptotic partition formulae, I:
Plane partitions}, Quart. J. Math. Oxford Ser. 2(1931), 177-189.

\end{thebibliography}
\end{document}